\documentclass[reqno,oneside,12pt,draft]{amsart}%
\usepackage{amsfonts}
\usepackage{amsmath}
\usepackage{amssymb}
\usepackage{graphicx}%
\hoffset - 2cm
\newtheorem{theorem}{Theorem}[section]
\theoremstyle{plain}

\newtheorem{lemma}{Lemma}[section]

\theoremstyle{definition}

\newtheorem{remark}{Remark}[section]

\numberwithin{equation}{section}

\begin{document}
\title[Lotka-Volterra systems ]{Existence of Periodic Solutions \\for the Lotka-Volterra type Systems }
\author{Norimichi Hirano}
\address{Department of Mathematics, Graduate School of Environment and Information
Sciences, Yokohama National University, Tokiwadai, Hodogayaku, Yokohama, Japan}
\email{hirano@hiranolab.jks.ynu.ac.jp}
\author{S{\l }awomir Rybicki$^{\dag}$}
\address{Faculty of Mathematics and Computer Science, Nicolaus Copernicus University,
PL-87-100 Toru\'{n}, ul. Chopina 12/18, Poland}
\email{Slawomir.Rybicki@mat.uni.torun.pl}
\thanks{$^{\dag}$ Partially supported by the Ministry of Education and  Science (Poland), under grant 1 PO3A 009 27}
\date{\today}

\begin{abstract}
In this paper we prove the existence of non-stationary periodic solutions of delay Lotka-Volterra equations. In the
proofs   we use the   $S^1$-degree due to Dylawerski et al. \cite{[DGJM]}.

\end{abstract}
\maketitle

\section{Introduction}

The aim of this paper is to prove the existence of non-stationary periodic solutions of autonomous delay
differential equations of Lotka-Volterra type
\begin{equation}
\left\{
\begin{array}
[c]{ll}%
\dot u_{1}(t)= & u_{1}(t)(r_{1}-a_{11}u_{1}(t-\tau)-a_{12} u_{2}%
(t-\tau)-\ldots-a_{1n}u_{n}(t-\tau)),\\
\dot u_{2}(t)= & u_{2}(t)(r_{2}-a_{21}u_{1}(t-\tau)-a_{22} u_{2}%
(t-\tau)-\ldots-a_{2n}u_{n}(t-\tau)),\\
& \ldots\\
\dot u_{n}(t)= & u_{n}(t)(r_{n}-a_{n1}u_{1}(t-\tau)-a_{n2} u_{2}%
(t-\tau)-\ldots-a_{nn}u_{n}(t-\tau)),
\end{array}
\right.  \label{P0}%
\end{equation}
where $n\geq1,\tau>0,$ $r_{1},\ldots,r_{n} \in\mathbb{R}, a_{ij} \in
\mathbb{R},$ for $i,j=1,\ldots,n.$

It is known that a broad class of problems in mathematical biology, economics
and mechanics are described in the form above with initial conditions%
\begin{equation}
\left\{
\begin{array}
[c]{ll}%
u_{i}(s)=\varphi_{i}(s),\quad & s \in[-\tau,0],\quad\varphi_{i}(0)>0,\\
\varphi_{i}\in C([-\tau,0],\mathbb{R}), & \qquad i=1,2,\ldots,n.
\end{array}
\right.  \label{I}%
\end{equation}
In case $n=1,$ the problem (\ref{P0}) is known as delay logistic equation. The existence and multiplicity of
solutions of delay logistic equation has been investigated by many authors (cf. Goparlsamy \cite{[GOPA]} and Hale
\cite{[HALE1]} and references therein). To compare with the method employed here with that for delay logistic
equation, we illustrate the proof for the existence of periodic solutions of the logistic equation
\begin{equation}
\label{L}\dot{u}(t)=\alpha u(t)(1-u(t-\tau)),
\end{equation}
where $\alpha>0.$ For each initial function $\varphi\in C([-\tau,0],\mathbb{R})$ with $\varphi(0)>0,$ one can find a
solution $u(\varphi)$ of \eqref{L} with initial value $u(\varphi)(s)=\varphi(s),$ $s\in[-\tau,0].$ We put
\[
z(\varphi,\alpha)=\min\left\{  t>0:u(\varphi)(t)=0,\overset{\cdot}{u}%
(\varphi)(t)>0\right\}
\]
and $[A(\alpha)\varphi](t)=u(\varphi)(t-z(\varphi,\alpha)-\tau)$ for $t>0.$ Then one can see that each fixed point
$u$ of $A(\alpha)$ is a periodic solution of \eqref{L}. The existence of the non-stationary fixed points of
$A(\alpha)$ is proved by combination of the Hopf bifurcation theorem and fixed point theorems (cf. Section 11.4 of
Hale \cite{[HALE1]}). For $\tau=1,$ it is known that $\alpha=\pi/2$ is the bifurcation point of solutions of
\eqref{L} and for each $\alpha > \pi/2,$ problem \eqref{L} has a non-stationary periodic solution.

On the other hand, it is natural to ask if there are multiple solutions of \eqref{L} for   sufficiently large
$\tau.$ The multiple existence of periodic solution of \eqref{L} for  sufficiently large $\tau$ also follows from
the Hopf bifurcation. In general the methods employed for delay logistic equation are not valid for (\ref{P0}) with
$n > 1.$

In this paper we  work with  the space of periodic functions instead of considering the initial value problem and
make use of the $S^1$-degree, see \cite{[DGJM]}, to prove the multiplicity of solutions of problem (\ref{P0}).
Applications of the degree for equivariant maps to the study of periodic solutions of a van der Pol system one can
find in \cite{[BAFAKR]}, \cite{[HR]}.

To avoid unnecessary complexity, we restrict ourselves to the case   $n=2,$ that is, we consider the coupled
equations of the form
\begin{equation}
\label{PW}\left\{
\begin{array}
[c]{ll}%
\dot u(t)= & u(t)(r_{1}-a_{11}u(t-\tau)-a_{12}v(t-\tau)),\\
\dot v(t)= & v(t)(r_{2}-a_{21}u(t-\tau)-a_{22}v(t-\tau)).
\end{array}
\right.
\end{equation}
Our argument does not depend on any specific property of $n=2.$ That is why our result is valid for $n\geq2$ with
modifications of assumptions for the case that $n\geq2.$

We impose the following conditions on matrix $A=\left[
\begin{array}
[c]{cc}%
a_{11} & a_{12}\\
a_{21} & a_{22}%
\end{array}
\right] :$

\bigskip

\noindent(A0) $a_{11},a_{12},a_{21},a_{22}, b_{1},b_{2}>0,$ where $\displaystyle\left[
\begin{array}
[c]{c}%
b_{1}\\
b_{2}%
\end{array}
\right]  =A^{-1}\left[
\begin{array}
[c]{c}%
r_{1}\\
r_{2}%
\end{array}
\right]  ,$

\noindent(A1) $\langle Ax,x \rangle> 0$ for all $x \in\mathbb{R}^{2}
\setminus\{0\},$

\noindent(A2) a matrix $\left[
\begin{array}
[c]{cc}%
b_{1} a_{11} & b_{1} a_{12}\\
b_{2} a_{21} & b_{2} a_{22}%
\end{array}
\right] $ possesses two real eigenvalues $\mu_{1}, \mu_{2} > 0.$

\begin{remark} Notice that for an arbitrary $n \in \mathbb{N}$ assumptions (A0), (A1), (A2) can be reformulated in
the following way

\noindent (A0) $a_{ij}, b_i > 0$ for $1 \leq i,j \leq n,$

\noindent (A1) $\langle Ax,x \rangle> 0$ for all $x \in\mathbb{R}^{n} \setminus\{0\},$

\noindent (A2) a matrix $\mathrm{diag}(b_1,\ldots,b_n) \cdot A$ possesses only real, positive eigenvalues
$\mu_1,\ldots,\mu_p,$ where $p>1.$ Moreover,  algebraic multiplicity of every eigenvalue $\mu_i$ is equal to its
geometric multiplicity.

\end{remark}

We can now formulate the main result of this article.

\begin{theorem}
\label{t} Fix $\tau> 0$ such that $\min\left\{ \frac{2\pi}{\mu_{1}},
\frac{2\pi}{\mu_{2}}\right\}  < \tau< \infty.$ Assume that there are
$n_{1},n_{2} \in\mathbb{N} \cup\{0\}$ such that $n_{1} \neq n_{2}$ and for
$i=1,2,$
\begin{equation}
\label{taucon}\displaystyle \frac{\pi}{2}+2n_{i}\pi< \mu_{i} \tau< \frac{\pi
}{2}+2(n_{i}+1)\pi.
\end{equation}
Under the above assumptions there is at least one non-stationary $\tau $-periodic solution of \eqref{PW}.
\end{theorem}

After this introduction our paper is organized as follows.

For the convenience of the reader in Section \ref{s1d} we have repeated the relevant material from \cite{[DGJM]}
without proofs, thus making our exposition self-contained.

In Section \ref{fs} we have performed a functional setting for our problem. This section is of technical nature.
Namely, applying   transformation of functions and fixing the period we have obtained  a parameterized problem
\eqref{P} which is equivalent to the original problem \eqref{PW}. Next we have defined a Banach space $E$ which is
an infinite-dimensional representation of the group $S^1,$ an open $S^1$-invariant subset $\Theta_0 \subset E$ and
an $S^1$-equivariant compact operator $F : (E \times \mathbb{R}^+) \times [0,1] \rightarrow E,$ see formula
\eqref{fq}, such that solutions of equation $F(((x_1,x_2),\lambda),1)=(x_1,x_2)$ in $\Theta_0 \times \mathbb{R}^+$
are exactly periodic solutions of problem \eqref{P}.

In Section \ref{hs1m} we have defined an open, bounded $S^1$-invariant subset $\Omega_{\lambda_{1} ,\lambda_{2}}
\subset \Theta_0 \times \mathbb{R}^+ \subset E \times \mathbb{R}^{+}$ such that the homotopy $Q -F(\cdot,\theta),$
defined by \eqref{fq} does not vanish on $\partial \Omega_{\lambda_{1},\lambda_{2}}.$ This allow us to simplify
computations of the $S^1$-degree of $Q -F(\cdot,1)$ on $\Omega_{\lambda_{1} ,\lambda_{2}},$ see Lemma
\ref{homotopy1}.

In Section \ref{proof} we have proved Theorem \ref{t}.

\section{$\mathrm{S^1}$-degree}
\label{s1d}

In this section we have compiled some basic facts on the $S^1$-degree defined in \cite{[DGJM]}. Let $S^1=\{z
\in\mathbb{C}: \: \mid z\mid=1\}=\{e^{i\cdot\theta}:\theta\in[0,2\pi)\}$ be the group with an action given by the
multiplication of complex numbers. For any fixed $m\in\mathbb{N}$ we denote by $\mathbb{Z}_{m}$ a cyclic group of
order $m$ and define homomorphism $\rho_{m}:S^1\rightarrow GL(2,\mathbb{R})$ as follows
\[
\rho_{m}\left(  e^{i\cdot\theta}\right)  =\left[
\begin{array}
[c]{rr}%
\cos(m\theta) & -\sin(m\theta)\\
\sin(m\theta) & \cos(m\theta)
\end{array}
\right]  .
\]

Let $E$ be a Banach space which is an $S^1$-representation. We denote by $Q : E \times\mathbb{R} \rightarrow E$ the
projection. For each closed subgroup $H$ of $S^1$ and each $S^1$-invariant subset $\Omega\subset E,$ we denote by
$\Omega^{H}$ the subset of fixed points of the action of $H$ on $\Omega.$ For given $a \in E,$ $S^1_{a}=\left\{ s\in
S^1: s \cdot a =a\right\}  $ is called the isotropy group of $a$ and the set $S^1 \cdot a=\left\{ s \cdot a : s\in
S^1\right\}  $ is called the orbit of $a$. Denote by $\Gamma_{0}$ the free abelian group generated by $\mathbb{N}$
and let $\Gamma =\mathbb{Z}_{2}\oplus\Gamma_{0}.$ Then $\gamma\in\Gamma$ means $\gamma =\left\{ \gamma_{r}\right\}
,$ where $\gamma_{0}\in$ $\mathbb{Z}_{2}$ and $\gamma_{r}\in\mathbb{Z}$ for $r\in\mathbb{N}.$

Fix an open, bounded $S^1$-invariant subset $\Omega\subset E \times \mathbb{R}$ and continuous $S^1$-equivariant
compact mapping $\Phi: cl(\Omega) \rightarrow E$ such that $(Q + \Phi)(\partial\Omega) \subset E \setminus\{0\}.$ In
this situation the $S^1$-degree $\mathrm{Deg}(Q + \Phi,\Omega)=\left\{  \gamma_{r}\right\} \in\Gamma,$ where
$\gamma_{0}=$ deg$_{S^1}(Q+\Phi,\Omega)$ and $\gamma_{r}=$ deg$_{\mathbb{Z} _{r}}(Q+ \Phi,\Omega), r \in\mathbb{N},$
has been defined in \cite{[DGJM]}.

\begin{theorem}
[{\cite{[DGJM]}}]\label{dgjm} Let $E$ be a Banach space which is a representation of the group $S^1,$ $\Omega_{0},
\Omega_{1}, \Omega_{2} \subset\Omega$ be open bounded, $S^1$-invariant subsets of $E \times \mathbb{R}.$ Assume that
$\Phi: cl(\Omega) \rightarrow E$ is a compact $S^1$-equivariant mapping such that $(Q+\Phi)(\partial\Omega) \subset
E \setminus\{0\}.$ Then there exists a $\Gamma$-valued function $\mathrm{Deg}(Q+\Phi,\Omega)$ called the
$S^1$-degree, satisfying the following properties:
\begin{enumerate}

\item[(a)] if $\deg_{H}(Q+\Phi,\Omega)\neq0,$ then $(Q+\Phi)^{-1}(0)\cap
\Omega^{H} \neq\emptyset,$

\item[(b)] if $(Q+\Phi)^{-1}(0)\cap\Omega\subset\Omega_{0},$ then
$\mathrm{Deg}(Q+\Phi,\Omega)=\mathrm{Deg}(Q+\Phi,\Omega_{0}),$

\item[(c)] if $\Omega_{1} \cap\Omega_{2} = \emptyset$ and $(Q +\Phi)^{-1}(0)
\cap\Omega\subset\Omega_{1} \cup\Omega_{2},$ then
\[
\mathrm{Deg}(Q+\Phi,\Omega)=\mathrm{Deg}(Q+\Phi,\Omega_{1}) + \mathrm{Deg}%
(Q+\Phi,\Omega_{2}),
\]

\item[(d)] if $h:cl(\Omega) \times[0,1] \rightarrow E$ is an $S^1$-equivariant compact homotopy such that
\linebreak $(Q+h)(\partial\Omega \times[0,1])\subset E\setminus\left\{ 0\right\} ,$ then $\mathrm{Deg}
(Q+h_{0},\Omega)=\mathrm{Deg}(Q+h_{1},\Omega). $
\end{enumerate}
\end{theorem}

Let $E_{0}, E$ be real Banach spaces. We denote by $K(E_{0},E)$ the set of compact operators $B : E_{0} \rightarrow
E.$ Fix $B_{0} \in K(E,E).$ A real number $\mu$ is called a characteristic value of $B_{0}$ if $\dim\ker(I-\mu
B_{0}) >0.$ Suppose now that $A_{0}=I-B_{0}$ is an invertible operator. Denote by $\left\{  \mu_{1}<\mu_{2}
<\cdot\cdot\cdot<\mu_{p}\right\}  $ the set of all characteristic values of $B_{0}$ contained in $[0,1.]$ We set
$\mathrm{sgn} A_{0}=(-1)^{d},$ where $\displaystyle d=\sum_{i=1}^{p} \dim \ker(I-\mu_{i} B_{0}).$

Fix $B_{1} \in K(E \times\mathbb{R},E)$ and define $A=Q + B_{1} : E \times\mathbb{R} \rightarrow E.$ Assume that $A$
is surjective. Since $A$ is a Fredholm operator of index $1$ and $A$ is surjective, $\dim\ker A =1.$ Fix $v \in\ker
A \setminus\{0\}$ and define a linear functional $\xi: E \times\mathbb{R} \rightarrow E$ such that $\xi(v)=1.$
Finally define an operator $A^{\sim} : E \times\mathbb{R} \rightarrow E \times\mathbb{R}$ by $A^{\sim}(w)=(A w,
\xi(w))$ and $\mathrm{sgn}(A,v)=\mathrm{sgn}(A^{\sim}).$

Assume additionally that $f = Q + \Phi\in C^{1}(cl(\Omega),E).$ Suppose that $a\in\Omega$ is such that there is $k
\in\mathbb{N}$ such that $S^1_{a}=\mathbb{Z}_{k},$ $f^{-1}(0) \cap\Omega= S^1 \cdot a \approx S^1 \slash
\mathbb{Z}_{k}$ and $Df(a) : E \times\mathbb{R} \rightarrow E$ is surjective. Let $f^{\mathbb{Z}_{k}} :
E^{\mathbb{Z}_{k}} \times\mathbb{R} \rightarrow E^{\mathbb{Z}_{k}}$ denotes the restriction of $f.$ Then $D
f^{\mathbb{Z}_{k}}(a) :E^{\mathbb{Z}_{k}} \times\mathbb{R} \rightarrow E^{\mathbb{Z}_{k}}$ is also surjective. We
denote by $v$ the tangent vector to the orbit $S^1 \cdot a$ at $a.$ Notice that $v \in\ker Df^{\mathbb{Z}_{k} }(a).$

This theorem ensures  the nontriviality of the $S^1$-index of the non-degenerate $S^1$-orbit $S^1 \cdot a.$

\begin{theorem}
[{\cite{[DGJM]}}]\label{dgjm2} Under the above assumptions $\mathrm{deg}_{\mathbb{Z}_{k}}(f,\Omega)=
\mathrm{sgn}(Df^{\mathbb{Z}_{k}}(a),v).$ Moreover, $\mathrm{deg}_{\mathbb{Z}_{k'}}(f,\Omega)= 0$ for every $k' > k.$
\end{theorem}

\section{Functional setting}

\label{fs}

Throughout the rest of this article we assume that assumptions (A0), (A1) and (A2) are fulfilled. Moreover, we fix
$\tau>0$ satisfying assumptions of Theorem \ref{t}.

In this section we convert problem \eqref{PW} to an equivalent problem \eqref{P}. Next we define   spaces on which
we will work and define a homotopy $F$ of $S^1$-invariant compact mappings. The study of periodic solutions of
problem \eqref{P} is equivalent to the study of fixed $S^1$-orbits of the operator $F(\cdot,1).$

By a transformation of functions (cf. \cite{[HALE1]}), one can see that problem \eqref{PW} is equivalent to the
problem
\begin{equation}
\left\{
\begin{array}
[c]{ll}%
\dot{u}(t)= & -(a_{11}u(t-\tau)+a_{12}v(t-\tau))(b_{1}+u(t)),\\
\dot{v}(t)= & -(a_{21}u(t-\tau)+a_{22}v(t-\tau))(b_{2}+v(t)).
\end{array}
\right.  \label{PV}%
\end{equation}
Let $(u,v)\in C(\mathbb{R},\mathbb{R}) \times C(\mathbb{R},\mathbb{R})$ be a pair of periodic
functions with period $T.$ Then by putting $\displaystyle\lambda=\frac{T}%
{2\pi},$ $x_{1}(t)=u(\lambda t)$ and $x_{2}(t)=v(\lambda t)$ for
$t\in\mathbb{R},$ we have that $x=(x_{1},x_{2})$ is a $2\pi$-periodic solution
of problem
\begin{equation}
\left\{
\begin{array}
[c]{ll}%
\dot{x}_{1}(t)= & -\lambda(a_{11}x_{1}(t-\tau\slash\lambda)+a_{12}x_{2}%
(t-\tau\slash\lambda))(b_{1}+x_{1}(t)),\\
\dot{x}_{2}(t)= & -\lambda(a_{21}x_{1}(t-\tau\slash\lambda)+a_{22}x_{2}%
(t-\tau\slash\lambda))(b_{2}+x_{2}(t)).
\end{array}
\right.  \label{P}%
\end{equation}

We will study the existence of  $2\pi$-periodic solutions of (\ref{P}) for some $\lambda>0$ instead of looking for
periodic solutions of (\ref{PV}). We note that each function $u:[0,2\pi]\rightarrow\mathbb{R}$ with $u(0)=u(2\pi)$
is extended to a $2\pi$-periodic function on $\mathbb{R}.$ Therefore we identify a $2\pi$-periodic function $u$ on
$\mathbb{R}$ with a function on $[0,2\pi]$ with $u(0)=u(2\pi).$

Define a Banach space
\[
\hat{E}=\left\{  x\in C([0,2\pi],\mathbb{R}):\int_{0}^{2\pi}\dot{x}
(t)^{2}dt<\infty,x(0)=x(2\pi):\mathrm{and}:\int_{0}^{2\pi}x(t)dt=0\right\}  .
\]
with a norm $\left\Vert \cdot\right\Vert $   given by $\displaystyle\left\Vert x\right\Vert
^{2}=\int_{0}^{2\pi}\dot{x} (t)^{2}+x(t)^{2}dt,\text{ for }x\in\hat{E}.$ Moreover, we put $\left\vert u\right\vert
_{\infty}=\sup\left\{  \left\vert u(t)\right\vert :t\in \lbrack0,2\pi]\right\}  $ for $u\in\hat{E}.$ Put $E=\hat{E}
\times\hat{E}$ and define an action $\rho:S^1 \times E\rightarrow E$ of the group $S^1$ as follows
\begin{equation}
\rho(e^{i\phi},(x_{1}(t),x_{2}(t)))=(x_{1}(t+\phi),x_{2}(t+\phi))\mod 2\pi
.\label{action}%
\end{equation}

Define an open $S^1$-invariant subset $\Theta_{0}\subset E$ as follows
\[
\Theta_{0}=\left\{  (x_{1},x_{2})\in E:-b_{i}<x_{i}(t)\quad\text{for }%
t\in\lbrack0,2\pi],i=1,2\right\}  .
\]

Next we define a homotopy of $S^1$-equivariant mappings $F:(E\times \mathbb{R}^{+})\times\lbrack0,1]\rightarrow E$
associated to the problem \eqref{P} such that if $((x_{1},x_{2}),\lambda)\in\Theta_{0}\times\mathbb{R}%
^{+}$ satisfies $F((x_{1},x_{2}),\lambda,1)=(x_{1},x_{2}),$ then $((x_{1}, x_{2}),\lambda)$ is a solution of problem
(\ref{P}).

Let $\beta:E\rightarrow\lbrack0,1]$ be a continuous mapping and $\mathcal{N}%
:E\times\lbrack0,1]\rightarrow E$ be a mapping defined by
\[
\mathcal{N}\left(  (x_{1},x_{2}),\theta\right)  (t)=\left(
\begin{array}
[c]{c}%
-(a_{11}x_{1}(t-\tau\slash\lambda)+a_{12}x_{2}(t-\tau\slash\lambda
))(b_{1}+\theta x_{1}(t))\\
-(a_{21}x_{1}(t-\tau\slash\lambda)+a_{22}x_{2}(t-\tau\slash\lambda
))(b_{2}+\theta x_{2}(t))
\end{array}
\right)
\]
for $((x_{1},x_{2}),\theta)\in E\times\lbrack0,1].$ We put
\[
c_{1}((x_{1},x_{2}),\theta)=\left(
\begin{array}
[c]{c}%
c_{11}((x_{1},x_{2}),\theta)\\
c_{12}((x_{1},x_{2}),\theta)
\end{array}
\right)  =-\frac{\lambda}{2\pi}\int_{0}^{2\pi}\beta(x_{1},x_{2})\mathcal{N}%
((x_{1},x_{2}),\theta)(s)ds
\]
and%
\begin{align*}
c_{2}((x_{1},x_{2}),\theta) &  =\left(
\begin{array}
[c]{c}%
c_{21}((x_{1},x_{2}),\theta)\\
c_{22}((x_{1},x_{2}),\theta)
\end{array}
\right)  \\
&  =-\frac{\lambda}{2\pi}\int_{0}^{2\pi}\int_{0}^{t}\beta(x_{1},x_{2}%
)\mathcal{N}((x_{1},x_{2}),\theta)(s)dsdt-\pi c_{1}((x_{1},x_{2}),\theta).
\end{align*}
We define a mapping $F:\left(  E\times\mathbb{R}^{+}\right)  \times \lbrack0,1]\rightarrow C(\mathbb{R},\mathbb{R})
\times C(\mathbb{R},\mathbb{R})$ by
\begin{equation}\label{opef}
F(((x_{1},x_{2}),\lambda),\theta)(t)=
\end{equation}
\[
= -\lambda\int_{0}^{t}\beta(x_{1}
,x_{2})\mathcal{N}((x_{1},x_{2}),\theta)(s)ds-tc_{1}((x_{1},x_{2}),\theta)-c_{2}((x_{1},x_{2}),\theta).
\]

From the definitions of $c_{1}$ and $c_{2},$ we can see that for all $((x_{1},x_{2},\lambda),\theta)\in\left(
E\times\mathbb{R}^{+}\right) \times\lbrack0,1]$ $\displaystyle F(((x_{1},x_{2}),\lambda),\theta
)(0)=F(((x_{1},x_{2}),\lambda),\theta)(2\pi)$ and $\displaystyle\int_{0}%
^{2\pi}F(((x_{1},x_{2}),\lambda),\theta)(t)dt=0$ holds.

Summing up, $F(((x_{1},x_{2}),\lambda),\theta) \in E$ for all $(((x_{1} ,x_{2}),\lambda),\theta) \in\left( E
\times\mathbb{R}^{+}\right) \times[0,1].$ It is also easy to see that $F : \left( E \times\mathbb{R} ^{+}\right)
\times[0,1] \rightarrow E$ is an $S^1$-equivariant compact mapping.

From the definition of $F,$ we find that $(((x_{1},x_{2}),\lambda),\theta) \in\left( \Theta_{0}
\times\mathbb{R}^{+}\right)  \times[0,1]$ satisfies
\begin{equation}
\label{fq}F(((x_{1},x_{2}),\lambda),\theta)=Q((x_{1},x_{2}),\lambda)
\end{equation}
if and only if
\begin{equation}
\left\{
\begin{array}
[c]{ll}%
\overset{\cdot}{x}_{1}(t)= & -\lambda\beta(x_{1},x_{2})(a_{11}x_{1}
(t-\tau\slash \lambda)+a_{12}x_{2}(t-\tau\slash \lambda))(b_{1}+\theta
x_{1}(t)) -c_{11} ((x_{1},x_{2}),\theta),\\
\overset{\cdot}{x}_{2}(t)= & -\lambda\beta(x_{1},x_{2})(a_{21}x_{1}%
(t-\tau\slash \lambda)+a_{22}x_{2}(t-\tau\slash \lambda))(b_{2}+\theta
x_{2}(t)) -c_{21} ((x_{1},x_{2}),\theta).
\end{array}
\right.  \label{Fthe}%
\end{equation}

We claim that system \eqref{Fthe} is equivalent to%

\begin{equation}
\left\{
\begin{array}
[c]{ll}%
\overset{\cdot}{x}_{1}(t)= & -\lambda\beta(x_{1},x_{2})(a_{11}x_{1}
(t-\tau\slash \lambda)+a_{12}x_{2}(t-\tau\slash \lambda))(b_{1}+\theta
x_{1}(t)),\\
\overset{\cdot}{x}_{2}(t)= & -\lambda\beta(x_{1},x_{2})(a_{21}x_{1}%
(t-\tau\slash \lambda)+a_{22}x_{2}(t-\tau\slash \lambda))(b_{2}+\theta
x_{2}(t)).
\end{array}
\right.  \label{Fthe1}%
\end{equation}


What is left is to show that $c_{i1} ((x_{1},x_{2}),\theta)=0$ for $i=1,2.$

\noindent Fix $i \in\{1,2\}$ and notice that
\begin{align*}
\frac{d}{dt} \ln(b_{i} +x_{i}(t))=  &  \frac{\dot{x}_{i}(t)}{b_{i}+x_{i}(t)}\\
=  &  -\lambda\beta(x_{1},x_{2} )(a_{i1}x_{1}(t-\tau\slash \lambda
)+a_{i2}x_{2}(t-\tau\slash \lambda))-\frac{c_{i1} ((x_{1},x_{2}),\theta
)}{b_{i}+x_{i}(t)}.
\end{align*}

\noindent Thus

\noindent$\ln(b_{i}+x_{i}(t))-\ln(b_{i}+x_{i}(s))$

\noindent$\displaystyle =-\lambda\int_{s}^{t}\beta(x_{1},x_{2})(a_{i1}%
x_{1}(w-\tau\slash \lambda)+a_{i2}x_{2}(w-\tau\slash \lambda)) dw -\int
_{s}^{t}\frac{c_{i1}((x_{1},x_{2}),\theta)}{b_{i}+x_{i}(w)} dw$

\noindent for $s,t\in\mathbb{R}$ with $s<t.$

Since $x_{i}(2\pi)=x_{i}(0), \displaystyle \int_{0}^{2\pi}x_{1}(t) dt
=\int_{0}^{2\pi}x_{2}(t) dt=0,$ $\displaystyle \int_{0}^{2\pi}\frac
{c_{i1}((x_{1},x_{2}),\theta)}{b_{i}+x_{i}(t)} dt=0. $

Finally, condition $x_{i}(t)>-b_{i}$ for all $t \in[0,2\pi],$ implies
$c_{i1}((x_{1},x_{2}),\theta)=0,$ which completes the proof.

\bigskip

We finish this section with the following lemma which yields apriori estimtes
for periodic solutions of problem \eqref{PA}.

\begin{lemma}
\label{MM}$(1)$ For $\lambda_{1},\lambda_{2}\in\mathbb{R}^{+}$ with
$\lambda_{1}<\lambda_{2},$ there exist positive numbers $m_{0},$ $\left\{
d_{i}\right\}  _{1\leq i\leq4}$ such that for each $\lambda\in[ \lambda
_{1},\lambda_{2}],\alpha\in[0,1]$ and $\tau\geq1,$ each solution $(x_{1}%
,x_{2})\in\Theta_{0}$ of the following problem
\begin{equation}
\left\{
\begin{array}
[c]{ll}%
\overset{\cdot}{x}_{1}(t)= & -\alpha\lambda(a_{11}x_{1}(t-\tau/\lambda
)+a_{12}x_{2}(t-\tau/\lambda))(b_{1}+x_{1}(t)),\\
\overset{\cdot}{x}_{2}(t)= & -\alpha\lambda(a_{21}x_{1}(t-\tau/\lambda
)+a_{22}x_{2}(t-\tau/\lambda))(b_{2}+x_{2}(t)),
\end{array}
\right.  \label{PA}%
\end{equation}
satisifes $\left\vert \overset{\cdot}{x}_{i}\right\vert _{\infty}<m_{0},$
$\left\vert \overset{\cdot\cdot}{x}_{i}\right\vert _{\infty}<m_{0}$ for
$i=1,2,$ and
\[
-b_{1} < -d_{1}<x_{1}(t)<d_{3}, -b_{2} < -d_{2}<x_{2}(t)<d_{4}\quad\text{on
}[0,2\pi];
\]

\noindent$(2)$ there exists $\alpha_{0}\in(0,1)$ such that there is no
nontrivial solution of $($\ref{PA}$)$ for $\alpha\in[0,\alpha_{0}]$
\end{lemma}

\begin{proof}
(1) Let $\lambda_{1},\lambda_{2}\in\mathbb{R}^{+}$ with $\lambda_{1}%
<\lambda_{2}$ and $\lambda\in\lbrack\lambda_{1},\lambda_{2}].$ Let
$(x_{1},x_{2})\in\Theta_{0}$ be a solution of (\ref{PA}). Then for
$i\in\{1,2\}$ we find that
\[
\ln(b_{i}+x_{i}(t))-\ln(b_{i}+x_{i}(s))=-\alpha\lambda\int_{s}^{t}\sum
_{j=1}^{2}a_{ij}x_{j}(w-\tau/\lambda)dw
\]
for $t,s\in\mathbb{R}$ with $s\leq t.$ Let $s\in\mathbb{R}$ such that
$x_{i}(s)=0.$ Then we have
\[
b_{i}+x_{i}(t)=b_{i}\exp(-\alpha\lambda\int_{s}^{t}\sum_{j=1}^{2}a_{ij}%
x_{j}(w-\tau/\lambda)dw).
\]
Then noting that $x_{i}(t)>-b_{i}$ for $i=1,2,$ and taking into account
assumption (A0) we obtain
\begin{equation}
x_{1}(t)<d_{3}=b_{1}\exp(2\pi\lambda(a_{11}b_{1}+a_{12}b_{2}))-b_{1}%
\qquad\text{for all }t\in\lbrack0,2\pi].\label{x1m}%
\end{equation}
Similarly as above we obtain
\begin{equation}
x_{2}(t)<d_{4}=b_{2}\exp(2\pi\lambda(a_{21}b_{1}+a_{22}b_{2}))-b_{2}%
\qquad\text{for all }t\in\lbrack0,2\pi].\label{x2m}%
\end{equation}

On the other hand, for all $t \in[0,2\pi]$ we have
\[
-b_{1} < -d_{1} = b_{1}\exp(-2\pi\lambda(a_{11}d_{3} + a_{12} d_{4}%
))-b_{1}<x_{1}(t)
\]
and
\[
-b_{2} < -d_{2}=b_{2}\exp(-2\pi\lambda(a_{21}d_{3}+a_{22}d_{4}))-b_{2} <
x_{2}(t).
\]
From (\ref{PA}) and the inequalities above we find that
\begin{equation}
\left\vert \dot{x}_{i}(t)\right\vert \leq\max_{i=1,2}C(\left\vert
a_{i1}\right\vert d_{3}+\left\vert a_{i2}\right\vert d_{4})(b_{i}+d_{i+2}
)\quad\text{for } t \in[0,1] \text{ and } i=1,2.\label{deri}%
\end{equation}
\qquad We also have by differentiating the both sides of (\ref{PA}) and using
the inequalities above that $\left\{  \left\vert \ddot{x}_{i}(t)\right\vert :t
\in[0,2\pi], i=1,2\right\} $ is bounded, which completes the proof of (1).

(2) Suppose that there exists a sequence $\left\{  \alpha_{n}\right\} \subset\mathbb{R}^{+},$ $\left\{
(x_{1n},x_{2n})\right\}  \subset E$ such that $\displaystyle \lim_{n\rightarrow\infty}\alpha_{n}=0$ and each
$(x_{1n},x_{2n})$ is a solution of (\ref{PA}) with $\alpha=\alpha_{n}.$ Then (\ref{x1m}) and (\ref{x2m}) holds with
$x_{1}$ and $x_{2}$ replaced by $x_{1n}$ and $x_{2n},$ respectively. Then we have that $\displaystyle
\lim_{n\rightarrow\infty}\left\vert x_{in}\right\vert _{\infty}=0$, $i=1,2.$ By subtracting subsequences, we may
assume, without any loss of generality, that $\left\vert x_{1n}\right\vert _{\infty}\geq\left\vert x_{2n}\right\vert
_{\infty}$ for all $n\geq1.$ We put $u_{in}(t)=x_{in}/\left\vert x_{1n}\right\vert _{\infty}$ for $i=1,2.$ Then we
have that $u_{in}\in H$ for $n\geq1$ and $i=1,2.$ We also have $\left\vert u_{1n}\right\vert _{\infty}=1$ for all
$n\geq1.$ Then it follows that for $n$ sufficiently large $\displaystyle \left\vert \dot{u}_{1n}(t)\right\vert
\leq2b_{1}\alpha _{n}\lambda(a_{11}+a_{12}) \quad\text{for all } t \in[0,2\pi] .$ That is $\displaystyle
\lim_{n\rightarrow\infty}\left\vert \dot{u}_{1n}\right\vert _{\infty}=0.$ This contradicts to the fact that
$\left\vert u_{1n}\right\vert _{\infty}=1$ for all $n \geq1,$ which completes the proof of (2).
\end{proof}

\section{Homotopies of admissible $S^1$-equivariant mappings} \label{hs1m}

The aim of this section is to define an open, bounded $S^1$-invariant subset $\Omega_{\lambda_{1} ,\lambda_{2}}
\subset \Theta_0 \times \mathbb{R}^+ \subset E \times \mathbb{R}^{+}$ such that the homotopy $Q -F(\cdot,\theta),$
defined by \eqref{fq}, does not vanish on $\partial \Omega_{\lambda_{1},\lambda_{2}}.$ We underline that solutions
of equation $Q((x_1,x_2),\lambda) =F((x_1,x_2),\lambda),1)$ in  $\Omega_{\lambda_{1} ,\lambda_{2}}$ are exactly the
periodic solutions of problem \eqref{P} in $\Omega_{\lambda_{1} ,\lambda_{2}}$.

We finish this section with Lemma \ref{homotopy1}, where we reduce the computation of the $S^1$-degree of $Q
-F(\cdot,1)$ on $\Omega_{\lambda_{1} ,\lambda_{2}}$ to the computation of  the $S^1$-degree of $Q -F(\cdot,0)$ on
$\Omega_{\lambda_{1} ,\lambda_{2}}$.

We first consider the following eigenvalue problem associated with problem (\ref{P})

\begin{equation}
\left\{
\begin{array}
[c]{ll}%
\overset{\cdot}{u}(t)= & -\gamma\lambda b_{1}(a_{11}u(t-\tau/\lambda
)+a_{12}v(t-\tau/\lambda)),\\
\overset{\cdot}{v}(t)= & -\gamma\lambda b_{2}(a_{21}u(t-\tau/\lambda
)+a_{22}v(t-\tau/\lambda)),
\end{array}
\right.  \label{linear}%
\end{equation}
where $\lambda,\tau>0,$ $\gamma\in\mathbb{R}$ and $(u,v)\in E.$  By assumption (\ref{taucon}) $\mu _{1}\neq\mu_{2}.$
Hence from assumption (A2) it follows that there exists a non-degenerate matrix $P$ such that
\begin{equation}\label{PAP}%
P\left[
\begin{array}
[c]{cc}%
a_{11}b_{1} & a_{12}b_{1}\\
a_{21}b_{2} & a_{22}b_{2}%
\end{array}
\right]  P^{-1}=\left[
\begin{array}
[c]{cc}%
\mu_{1} & 0\\
0 & \mu_{2}%
\end{array}
\right].
\end{equation}
 Then linear problem (\ref{linear}) is transformed
into the form
\begin{equation}
\left\{
\begin{array}
[c]{lll}%
\overset{\cdot}{u_{1}}(t)= & -\gamma\lambda\mu_{1}u_{1}(t-\tau/\lambda), & \\
\overset{\cdot}{u_{2}}(t)= & -\gamma\lambda\mu_{2}u_{2}(t-\tau/\lambda)), &
\end{array}
\right.  \label{linearv}%
\end{equation}
by putting $(u_{1},u_{2})=P^{-1}(u,v).$

Let $i=1,2.$ We put $\displaystyle u_{i}%
(t)=\sum_{k=1}^{\infty}(c_{k}\cos kt+s_{k}\sin kt),$ where $\left\{ c_{k}\right\}  ,\left\{  s_{k}\right\}
\subset\mathbb{R}.$ Then \linebreak $\displaystyle\dot{u}_{i}(t)=\sum_{k=1}^{\infty}\left(  (ks_{k})\cos
kt+(-kc_{k})\sin kt\right)  $ and
\begin{align*}
-\gamma\lambda\mu_{i}u(t-\tau^{\prime}) &  =-\gamma\lambda\mu_{i}\sum
_{k=1}^{\infty}\{c_{k}(\cos kt\cos k\tau^{\prime}+\sin k\tau^{\prime}\sin
kt)\\
&  \qquad\qquad\qquad+s_{k}(\cos k\tau^{\prime}\sin kt-\sin k\tau^{\prime}\cos
kt)\}\\
&  =-\gamma\lambda\mu_{i}\sum_{k=1}^{\infty}\{(c_{k}\cos k\tau^{\prime}%
-s_{k}\sin k\tau^{\prime})\cos kt\\
&  \qquad\qquad\qquad+(c_{k}\sin k\tau^{\prime}+s_{k}\cos k\tau^{\prime})\sin
kt\},
\end{align*}
where $\tau^{\prime}=\tau/\lambda.$ If $u_1$ is a nontrivial solution of (\ref{linearv}), then
\begin{align*}
ks_{k} &  =-\gamma\lambda\mu_{1}(c_{k}\cos k\tau^{\prime}-s_{k}\sin
k\tau^{\prime}),\\
-kc_{k} &  =-\gamma\lambda\mu_{1}(c_{k}\sin k\tau^{\prime}+s_{k}\cos k\tau^{\prime}),
\end{align*}
for all $k\in\mathbb{N}.$ That is we obtain the following system of linear
equations
\begin{align*}
c_{k}\cos k\tau^{\prime}-s_{k}\left(  \sin k\tau^{\prime}-\frac{k}%
{\gamma\lambda\mu_{1}}\right)   &  =0,\\
c_{k}\left(  \sin k\tau^{\prime}-\frac{k}{\gamma\lambda\mu_{1}}\right) +s_{k}\cos k\tau^{\prime} &  =0,
\end{align*}
for all $k\in\mathbb{N}.$ Then $\cos^{2}k\tau^{\prime}+\left(  \sin
k\tau^{\prime}-\frac{k}{\gamma\lambda\mu_{1}}\right)  ^{2}=0$ and therefore we
find that%
\begin{equation}
\frac{k}{\gamma\lambda\mu_{1}}=1,\qquad\frac{k\tau}{\lambda}=\frac{\pi}%
{2}+2n\pi,\quad\text{for some }n\in\mathbb{N}\cup\{0\}.\label{eigen1}%
\end{equation}
If $u_2$ is a nontrivial solution of (\ref{linearv}), then by the same argument as above we obtain that
\begin{equation}
\frac{k}{\gamma\lambda\mu_{2}}=1,\qquad\frac{k\tau}{\lambda}=\frac{\pi}%
{2}+2n\pi,\quad\text{for some }n\in\mathbb{N}\cup\{0\}.\label{eigen2}%
\end{equation}
Consequently, we have that the eigenvalue $\gamma$ of problem (\ref{linearv}) is of the form
\begin{equation}
\gamma=\frac{1}{\mu_{i}\tau}\left(  \frac{\pi}{2}+2n\pi\right)  ,\quad i=1,2,\text{ and
}n\in\mathbb{N}.\label{eigenvalue}
\end{equation}
Based on the observation above we obtain the following lemma.

\begin{lemma}
\label{mm}Let $\lambda_{1},\lambda_{2}\in\mathbb{R}^{+}$ with $\lambda
_{1}<\lambda_{2}.$ Then the origin is an isolated solution of \eqref{P} i.e.
there exists $m_{1} > 0$ such that if $((x_{1},x_{2}),\lambda)\in\Theta_{0}
\times[\lambda_{1},\lambda_{2}]$ is a nontrivial solution of \eqref{P} then
$(x_{1},x_{2}) \notin\left\{ (x_{1},x_{2})\in\Theta_{0} : \left\Vert
x_{1}\right\Vert \leq m_{1},\left\Vert x_{2}\right\Vert \leq m_{1}\right\}  .
$
\end{lemma}

\begin{proof}
Let $\lambda_{1},\lambda_{2}\in\mathbb{R}^{+}$ with $\lambda_{1}<\lambda_{2}.$
Suppose that there exists a sequence \linebreak$\left\{  ((x_{1n},x_{2n})
,\lambda_{n})\right\}  \subset E\times\mathbb{R}^{+}$ such that each
$((x_{1n},x_{2n}),\lambda_{n})$ is a solution of (\ref{P}) and
\[
\lim_{n\rightarrow\infty}\left\Vert x_{1n}\right\Vert =\lim_{n\rightarrow
\infty}\left\Vert x_{2n}\right\Vert =0.
\]
We may assume that $\displaystyle \lim_{n\rightarrow\infty}\lambda_{n}%
=\lambda_{0} \in[\lambda_{1},\lambda_{2}]$ and $\left\vert x_{2n}\right\vert
_{\infty}\leq\left\vert x_{1n}\right\vert _{\infty}$ for all $n\geq1.$ We put
$u_{in}(t)=x_{in}/\left\vert x_{1n}\right\vert _{\infty}$ for each $n\geq1$
and $i=1,2.$ Then we have
\[
\left\{
\begin{array}
[c]{ll}%
\dot{u}_{1n}(t)= & -\lambda_{n}(a_{11}u_{1n} (t-\tau/\lambda_{n})+a_{12}%
u_{2n}(t-\tau/\lambda_{n}))(b_{1} +x_{1n}(t)),\\
\dot{u}_{2n}(t)= & -\lambda_{n}(a_{21}u_{1n} (t-\tau/\lambda_{n})+a_{22}%
u_{2n}(t-\tau/\lambda_{n}))(b_{2} +x_{2n}(t)).
\end{array}
\right.
\]
Then one can see that $\sup\left\{  \left\Vert \dot{u}_{in}\right\Vert : n
\geq1, i=1,2\right\}  < \infty.$ By differentiating the equalities above, we
also have that $\sup\left\{  \left\Vert \ddot{u}_{in}\right\Vert : n \geq1,
i=1,2\right\}  <\infty.$ Therefore we may assume that for each $i,$ $u_{in}
\rightarrow u_{i}$ and $\dot{u}_{in} \rightarrow\dot{u}_{i}$ strongly in $E.$
Then we have
\begin{equation}
\left\{
\begin{array}
[c]{ll}%
\dot{u}_{1}(t)= & -\lambda_{0} b_{1}(a_{11}u_{1}(t-\tau/\lambda_{0}
)+a_{12}u_{2}(t-\tau/\lambda_{0})),\\
\dot{u}_{2}(t)= & -\lambda_{0} b_{2}(a_{21}u_{1}(t-\tau/\lambda_{0}
)+a_{22}u_{2}(t-\tau/\lambda_{0})).
\end{array}
\right.  \label{linear2}%
\end{equation}
That is (\ref{linear}) holds with $\gamma=1.$ By the assumption
\eqref{taucon}, we have that \eqref{eigenvalue} does not hold with $\gamma=1.$
Therefore problem \eqref{linear2} has no nontrivial solution. Then
$u_{1}\equiv0.$ This contradicts the definition of $u_{1}.$
\end{proof}

Now fix $\lambda_{1},\lambda_{2}\in\mathbb{R}^{+}$ with $\lambda_{1} <\lambda_{2}$ and $\alpha_{0}, m_{0},$ $m_{1},$
$\left\{  d_{i}\right\} _{1\leq i\leq4}$ be the positive numbers satisfying the assertion of Lemma \ref{MM} and
Lemma \ref{mm}. We may assume without any loss of generality that $d_{i+2} >(b_{i}+d_{i})/2$ for $i=1,2.$ Define
open bounded $S^1$-invariant subsets
\[
\Theta_{M}=\left\{  (x_{1},x_{2})\in E:-\frac{b_{i}+d_{i}}{2}<x_{i}
(t)<2d_{i+2},\quad\text{for }i=1,2,t\in[0,2\pi]\right\} ,
\]
\[
\widetilde{\Theta}_{M}=\left\{  (x_{1},x_{2})\in E:-d_{i}<x_{i}(t)<d_{i+2},
\quad\text{for }i=1,2,t\in[0,2\pi]\right\} ,
\]
and closed $S^1$-invariant subset as follows
\[
\Theta_{m_{1}}=\left\{  (x_{1},x_{2})\in E:\left\Vert x_{i}\right\Vert \leq
m_{1}/2,\quad i=1,2\right\} .
\]
Since $b_{i} > d_{i}, i=1,2$ and $m_{1} > 0$ can be chosen sufficiently small
\[
\Theta_{m_{1}} \subsetneq\widetilde{\Theta}_{M} \subsetneq\Theta_{M}
\subsetneq\Theta_{0}.
\]
Moreover, define open bounded $S^1$-invariant subsets in the following way
\[
\Theta_{1}=\left\{  (x_{1},x_{2})\in E:\int_{0}^{2\pi} \dot{x}_{i}(t)^{2} dt <2\pi m_{0}^{2},\quad i=1,2\right\}
,\widetilde{\Theta} = (\widetilde{\Theta}_{M} \cap\Theta_{1}) \setminus \Theta_{m_{1}},
\]
\begin{equation}\label{Theta}
\Theta=(\Theta_{M} \cap\Theta_{1}) \setminus\Theta_{m_{1}},
\end{equation}
and notice that $\widetilde{\Theta}_{M} \subset cl(\widetilde{\Theta}_{M})
\subset\Theta.$ Then we can choose $\delta_{0}>0$ such that $\mathrm{dist}%
^{2}(\widetilde{\Theta}_{M},\partial\Theta) \geq\delta_{0}.$  Let $\xi:[0,+
\infty) \rightarrow(0,1]$ be a smooth function such that
\begin{equation}
\xi(t)=\left\{
\begin{array}
[c]{lcl}%
\alpha_{0} &  & \text{for }t=0,\\
\text{strictly increasing } &  & \text{for } 0 < t < \delta_{0},\\
1 &  & \text{for }t\geq\delta_{0}.
\end{array}
\right. \label{kusai}%
\end{equation}
We put that%
\begin{equation}
\beta(x_{1},x_{2})=\xi\left( \mathrm{dist}^{2}((x_{1},x_{2}),\partial
\Theta)\right)  \qquad\text{for }(x_{1},x_{2})\in E.\label{betadef}%
\end{equation}
Then $\beta\in C^{1}(E ;\mathbb{R)}$ and we have
\begin{equation}
\beta(x_{1},x_{2})=\left\{
\begin{array}
[c]{ccl}%
1 &  & \text{for }(x_{1},x_{2})\in\Theta_{m_{1}},\\
\alpha_{0} &  & \text{for }(x_{1},x_{2}) \in cl(E \setminus(\Theta_{M}
\cap\Theta_{1})).
\end{array}
\right.  \label{beta}%
\end{equation}
Put $\Omega_{\lambda_{1},\lambda_{2}}=\Theta\times(\lambda_{1},\lambda_{2})$
and
\[
\mathcal{S}=\left\{  (((x_{1},x_{2}),\lambda),\theta) \in\Omega_{\lambda
_{1},\lambda_{2}} \times[0,1] : F(((x_{1},x_{2}),\lambda),\theta
)=Q((x_{1},x_{2}),\lambda) \right\} .
\]

\begin{lemma} \label{mm2}
Under the above assumptions \label{S} $\mathcal{S} \cap((\partial\Theta
\times(\lambda_{1},\lambda_{2})) \times[0,1]) = \emptyset.$
\end{lemma}

\begin{proof}
Suppose that $(((x_{1},x_{2}),\lambda),\theta) \in\mathcal{S}.$ Then
(\ref{Fthe}) holds for $(((x_{1},x_{2}),\lambda),\theta).$ Multiplying
(\ref{Fthe}) by $\theta,$ and denoting $u_{i}(t)=\theta x_{i}(t)$  for
$i=1,2,$ we find that
\[
\left\{
\begin{array}
[c]{ll}%
\dot u_{1}(t)= & -\lambda\beta(x_{1},x_{2})(a_{11}u_{1} (t-\tau/\lambda) +
a_{12}u_{2}(t-\tau/\lambda))(b_{1}+u_{1}(t)),\\
\dot u_{2}(t)= & -\lambda\beta(x_{1},x_{2})(a_{21}u_{1} (t-\tau/\lambda
)+a_{22}u_{2}(t-\tau/\lambda)((b_{2}+u_{2}(t)),
\end{array}
\right.
\]
holds. By Lemma \ref{MM} we obtain $\displaystyle \int_{0}^{2\pi} \dot{x}%
_{i}(t)^{2} dt <2\pi m_{0}^{2},$ $i=1,2$. Then we have $(x_{1},x_{2}%
)\notin\partial\Theta_{1}.$

If $(x_{1},x_{2})\in\partial\Theta_{M},$ then we have that $\beta(x_{1}%
,x_{2})=\alpha_{0}.$ Then by (2) of Lemma \ref{MM}, we have that $u_{1}\equiv
u_{2}\equiv0.$ This contradicts to $(x_{1},x_{2})\in cl(\Theta).$ If
$(x_{1},x_{2})\in\partial\Theta_{m_{1}},$ then $\beta(x_{1},x_{2})=1$ and
$\left\Vert u_{i}\right\Vert <m_{1}/2,$ for $i=1,2.$ Then by Lemma \ref{mm},
we have that $x_{1}\equiv x_{2}\equiv0.$ Thus we have that $\mathcal{S}
\cap\partial\Omega_{\lambda_{1},\lambda_{2}}=\emptyset,$ which completes the proof.
\end{proof}

\noindent The following result is known. For completeness, we give a proof.

\begin{lemma}
\label{okure}Suppose that $\tau=0.$ Then problem \eqref{Fthe} does not have
non-stationary periodic solution $((x_{1},x_{2}),\lambda)\in E \times
\mathbb{\mathbb{R}^{+}}$ for any $\theta\in[0,1].$
\end{lemma}

\begin{proof}
Let $\theta\in\lbrack0,1]$ and $((x_{1},x_{2}),\lambda)\in\Theta_{0}%
\times\mathbb{R}^{+}$ satisfy (\ref{Fthe}). We first consider the case that
$\theta>0.$ Since $\tau=0,$ problem (\ref{Fthe}) reduces to the problem
\begin{equation}
\left\{
\begin{array}
[c]{ll}%
\overset{\cdot}{x}_{1}(t)= & -\lambda(a_{11}x_{1}(t)+a_{12}x_{2}%
(t))(b_{1}+\theta x_{1}(t)),\\
\overset{\cdot}{x}_{2}(t)= & -\lambda(a_{21}x_{1}(t)+a_{22}x_{2}%
(t))(b_{2}+\theta x_{2}(t)).
\end{array}
\right.  \label{without}%
\end{equation}
We integrate the both sides of \eqref{without} from $0$ to $2\pi.$ Then by the
periodicity, we have
\[
\left\{
\begin{array}
[c]{l}%
\displaystyle a_{11}\int_{0}^{2\pi}x_{1}(t)^{2}dt+a_{12}\int_{0}^{2\pi}%
x_{1}(t)x_{2}(t)dt=0,\\
\displaystyle a_{21}\int_{0}^{2\pi}x_{1}(t)x_{2}(t)dt+a_{22}\int_{0}^{2\pi
}x_{2}(t)^{2}dt=0.
\end{array}
\right.
\]
Then one can see that $x_{1}=x_{2}\equiv0$ from the condition (A1). We next
consider the case that $\theta=0.$ In this case we multiply equations
(\ref{without}) by $x_{i}$ and integrate over $[0,2\pi].$ Then we have the
equalities above. This completes the proof.
\end{proof}

\begin{lemma}
\label{homotopy1} Suppose that $\displaystyle \lambda_{1}=\frac{\tau} {2j_{1}\pi} < \lambda_{2}=\frac{\tau}{2
j_{2}\pi},$ where $j_{1}, j_{2} \in\mathbb{N}.$ Then
\[
\mathrm{Deg}(Q-F(\cdot,0),\Omega_{\lambda_{1},\lambda_{2}})=\mathrm{Deg} (Q-F(\cdot,1)
,\Omega_{\lambda_{1},\lambda_{2}}).
\]

\end{lemma}

\begin{proof}
To prove the assertion, it is sufficient to show that there exists no solution of (\ref{Fthe}) in
$\partial\Omega_{\lambda_{1},\lambda_{2}}=\partial\left( \Theta\times(\lambda_{1},\lambda_{2})\right) =cl(\Theta)
\times\{\lambda _{1},\lambda_{2}\} \cup\partial\Theta\times(\lambda_{1},\lambda_{2}).$ We first see that there
exists no solution on $cl(\Theta) \times\{\lambda _{1},\lambda_{2}\}.$ From the definitions of
$\lambda_{1},\lambda_{2},$ we have that the problem (\ref{Fthe}) is equivalent to (\ref{without}) with
$\lambda=\lambda_{1}$ or $\lambda=\lambda_{2}.$ Then by Lemma \ref{okure}, we find that $x_{1}=x_{2}=0.$ This
contradicts to the assumption that $(x_{1},x_{2})\in cl(\Theta).$ We also have by Lemma \ref{S} that there exists no
solution of (\ref{Fthe}) in $\partial\Theta\times(\lambda_{1},\lambda _{2}),$ which completes the proof.
\end{proof}



\section{Proof of Theorem \ref{t}} \label{proof}

Throughout this section we assume that $\displaystyle\lambda_{1}=\frac{\tau }{2j_{1}\pi}<\lambda_{2}=
\frac{\tau}{2j_{2}\pi},$ where $j_{1},j_{2} \in\mathbb{N}$ and put $\Omega_{\lambda_{1},\lambda_{2}}=\Theta\times
(\lambda_{1},\lambda_{2}).$ From Theorem \ref{dgjm} and Lemma 4.4 it follows that to finish the proof of Theorem
1.1, it is sufficient to show that $\mathrm{Deg}(Q-F(\cdot,0),\Omega_{\lambda_{1},\lambda_{2}})\neq \Theta \in
\Gamma.$ But the mapping $F(\cdot,0)$ is still too complicated to calculate the $S^1$-degree. Therefore we will
provide another homotopy $G$ of $S^1$-equivariant compact mappings such that $F(\cdot ,0)=G(\cdot,0)$ and the
$S^1$-degree of of $Q-G(\cdot,1)$ on $\Omega_{\lambda_{1},\lambda_{2}} $ is easy to compute.

We fix a $C^{1}$-mapping $\sigma:[\lambda_{1} ,\lambda_{2}]\rightarrow\lbrack\lambda_{1},\lambda_{2}]$ such that
$\sigma$ is increasing on $[\lambda_{1},\lambda_{2}]$ with $\sigma(\lambda_{1})=\lambda_{1}$ and
$\sigma(\lambda_{2})=\lambda_{2},$ and
\[
\sigma(\lambda_{k,n})=\lambda_{k,n}\text{ and }\dot{\sigma}(\lambda
_{k,n})=0\qquad\text{for each }\lambda_{k,n}=\frac{k\tau}{\frac{\pi}{2}+2n\pi
}\in\lbrack\lambda_{1},\lambda_{2}],\quad k,n\in\mathbb{N}.
\]
We now define a homotopy of $S^1$-equivariant mappings
$G:\Omega_{\lambda_{1},\lambda_{2}}\times\lbrack0,1]\rightarrow E$ by
\begin{equation}
G(((x_{1},x_{2}),\lambda),\theta)=-(\theta\sigma(\lambda)+(1-\theta
)\lambda)\int_{0}^{t}\beta(x_{1},x_{2})\mathcal{N}((x_{1},x_{2}%
),0)ds.\label{G}%
\end{equation}
By definition of $\mathcal{N}(\cdot,0),$ we have $G(((x_{1},x_{2}%
),\lambda),\theta)\in E$ for $(((x_{1},x_{2}),\lambda),\theta)\in\left(
E\times\mathbb{R}^{+}\right)  \times\lbrack0,1].$ If $(((x_{1},x_{2}%
),\lambda),\theta)\in\Omega_{\lambda_{1},\lambda_{2}}\times\lbrack0,1]$
satisfies $Q((x_{1},x_{2}),\lambda)=G(((x_{1},x_{2}),\lambda),\theta)$ then
\begin{equation}
\left\{
\begin{array}
[c]{l}%
\dot{x}_{1}(t)=-(\theta\sigma(\lambda)+(1-\theta)\lambda)b_{1}\beta
(x_{1},x_{2})(a_{11}x_{1}(t-\tau/\lambda)+a_{12}x_{2}(t-\tau/\lambda)),\\
\dot{x}_{2}(t)=-(\theta\sigma(\lambda)+(1-\theta)\lambda)b_{2}\beta
(x_{1},x_{2})(a_{21}x_{1}(t-\tau/\lambda)+a_{22}x_{2}(t-\tau/\lambda)).
\end{array}
\right.  \label{gs}%
\end{equation}

\begin{lemma}
\label{homotopy2} Under the above assumptions:
\[
\mathrm{Deg}(Q-F(\cdot,0),\Omega_{\lambda_{1},\lambda_{2}})=\mathrm{Deg}
(Q-G(\cdot,1),\Omega_{\lambda_{1},\lambda_{2}}).
\]

\end{lemma}

\begin{proof}
By the same argument as in the proof of Lemma \ref{homotopy1}, we see that
\[
\mathrm{Deg}(Q-G(\cdot,0),\Omega_{\lambda_{1},\lambda_{2}})=\mathrm{Deg}%
(Q-G(\cdot,1) ,\Omega_{\lambda_{1},\lambda_{2}}).
\]
Then since $G(\cdot,0)=F(\cdot,0),$ we have by Lemma \ref{homotopy1} that the
assertion holds.
\end{proof}

For $n, m \in\mathbb{N}$ define $\displaystyle \Phi(n,m)=\left\{j \in\mathbb{N} : \left[ \frac{n}{j} \right]  <
\left[ \frac{m}{j} \right] = \frac{m}{j}\right\} .$ Notice that if $n < m$ then $\Phi(n,m) \neq\emptyset.$

The following lemma plays crucial role in our article.

\begin{lemma}
\label{homotopy3} Let assumptions of Theorem \ref{t} be fulfilled. If $n_{1} < n_{2}, j \in\Phi(n_{1},n_{2})$ and
$\displaystyle \lambda_{1}=\frac{\tau}{2(j+1)\pi}, \lambda_{2} =\frac{\tau}{2j\pi}$ then
$\mathrm{Deg}(Q-G(\cdot,1),\Omega_{\lambda_{1},\lambda_{2}}) \neq\Theta\in \Gamma.$
\end{lemma}
\begin{proof} Before we prove this lemma, we outline the main steps of the proof.
Namely, we  will prove that $(Q -G(\cdot,1))^{-1}(0) \cap \Omega_{\lambda_1 \lambda_2}$ consists of a finite number
of non-degenerate orbits $S^1 \cdot a_1,\ldots, S^1  \cdot a_p.$ Since these orbits are non-degenerate, there are
open bounded $S^1$-invariant subsets $U_i \subset cl(U_i) \subset \Omega_{\lambda_1 \lambda_2},i=1,\ldots,p,$ such
that $(Q -G(\cdot,1))^{-1}(0) \cap U_i=S^1  \cdot a_i, i=1, \ldots,p.$ Moreover, we will prove that there are $k_0
\in \mathbb{N}$ and $1 \leq i_0 \leq p$ such that $S^1_{a_{i_0}}=\mathbb{Z}_{k_0}$ and $S^1_{a_{i}} \not =
\mathbb{Z}_{k_0}$ for every $i \neq i_0.$

By Theorem \ref{dgjm} we obtain
$$\mathrm{Deg}(Q-G(\cdot,1),\Omega_{\lambda_{1},\lambda_{2}})=$$ $$=
\mathrm{Deg}(Q-G(\cdot,1),U_1)+ \ldots + \mathrm{Deg}(Q-G(\cdot,1),U_p) \in \Gamma.$$

From the above and Theorem \ref{dgjm2} we obtain that
$$\mathrm{Deg}_{\mathbb{Z}_{k_0}}(Q-G(\cdot,1),\Omega_{\lambda_{1},\lambda_{2}})=$$ $$=
\mathrm{Deg}_{\mathbb{Z}_{k_0}}(Q-G(\cdot,1),U_1)+ \ldots + \mathrm{Deg}_{\mathbb{Z}_{k_0}}(Q-G(\cdot,1),U_p) =
\mathrm{Deg}_{\mathbb{Z}_{k_0}}(Q-G(\cdot,1),U_{i_0})  \not = 0 \in \mathbb{Z}.$$

Let us begin the proof. First of all notice that since $\mu_{1} \neq\mu_{2}$ applying change of coordinates
\eqref{PAP} to the system \eqref{gs} we obtain the following equivalent system:
\begin{equation}
\left\{
\begin{array}
[c]{ll}%
\dot x_{1}(t)= & -\sigma(\lambda)\beta(P^{-1}(x_{1},x_{2})) \mu_{1}
x_{1}(t-\tau/\lambda),\\
\dot x_{2}(t)= & -\sigma(\lambda)\beta(P^{-1}(x_{1},x_{2})) \mu_{2}
x_{2}(t-\tau/\lambda)).
\end{array}
\right.  \label{Trans}%
\end{equation}

Notice that system \eqref{Trans} does not have solutions on $\partial
(P\Omega_{\lambda_{1},\lambda_{2}})=\partial(P\Theta\times(\lambda_{1}%
,\lambda_{2})).$ Since $\beta(P^{-1}(x_{1},x_{2})) = \alpha_{0}$ for any
$(x_{1},x_{2}) \in P(cl(E \setminus(\Theta_{M} \cap\Theta_{1}))),$ we find
that \eqref{Trans} does not have solutions on $P(E \setminus(\Theta_{M}
\cap\Theta_{1})) \times(\lambda_{1},\lambda_{2}).$

Therefore we can choose $R \gg r > 0$ such that
\[
\text{Deg}(Q-PG(\cdot,1)P^{-1},P\Theta\times(\lambda_{1},\lambda_{2}))=
\text{Deg} (Q-PG(\cdot,1)P^{-1},(D_{R}(E) \setminus cl(D_{r}(E))
\times(\lambda_{1},\lambda_{2})),
\]
where
\[
P(cl(\Theta_{M} \cap\Theta_{1})) \subset D_{R}(E)=\{x \in E : \|x\| < R\},
cl(D_{r}(E))=\{x \in E : \|x\| \leq r\} \subset P(\Theta_{m_{1}}).
\]

Here we replace $\beta$ by a function for which the calculation of degree is easier.

Let $\widetilde{\xi} \in C^{\infty}([0,+\infty),[0,1])$ with $\widetilde{\xi
}(t)=1$ for $0 \leq t \leq\sqrt{r},\widetilde{\xi}(t)=\alpha_{0}$ for $t
\in[\sqrt{R},+\infty)$ and $\widetilde{\xi}$ is strictly monotone decreasing
on $[\sqrt{r},\sqrt{R}].$ Define $\widetilde{\beta} : E \rightarrow[\alpha
_{0},1]$ as follows
\begin{equation}
\label{wer2}\widetilde{\beta}(x_{1},x_{2}):=\widetilde{\xi}(\|(x_{1}%
,x_{2})\|^{2}).
\end{equation}

We denote by $\widetilde{G}(\cdot,1)$ the mapping $G(\cdot,1)$ with $\beta$
replaced by $\widetilde{\beta}.$ Since maps $\beta, \widetilde{\beta}$
coincide on $\partial(D_{R}(E) \setminus cl(D_{r}(E)),$ $\widetilde{G}%
(\cdot,1)_{\mid\partial(D_{R}(E) \setminus cl(D_{r}(E))}=G(\cdot ,1)_{\mid\partial(D_{R}(E) \setminus cl(D_{r}(E))}$
and by the homotopy invariance of  the $S^1$-degree  we have
\[
\text{Deg}(Q-PG(\cdot,1)P^{-1},P\Theta\times(\lambda_{1},\lambda_{2}))=
\text{Deg} (Q- P\widetilde{G}(\cdot,1)P^{-1},(D_{R}(E) \setminus cl(D_{r}(E))
\times(\lambda_{1},\lambda_{2})).
\]
For the simplicity of notation we will denote $\widetilde{G}(\cdot,1),
\widetilde{\beta}$ and $\widetilde{\xi}$ by $G(\cdot,1), \beta$ and $\xi,$ respectively.

If $a(t)=(a_{1}(t),a_{2}(t))$ satisfies (\ref{Trans}), then we have by (\ref{eigenvalue}) that $\displaystyle
\lambda=\lambda_{k,n}=\frac{k\tau}{\frac{\pi}{2}+2n\pi}$ for some $k, n \in\mathbb{N}$ and $a_{i}(t)
\in\mathrm{span} \left\{  \cos kt,\sin kt\right\}  $ for $i=1,2.$ Then by the definition of $\sigma,$ we find that
$\displaystyle \sigma(\lambda_{k,n})=\lambda_{k,n}=\frac{k\tau}{\frac{\pi} {2}+2n\pi}$ for some $k,n\geq1.$

Now suppose that $a_{1} \not \equiv 0.$ Then since $\sigma(\lambda_{k,n} )\beta(a_{1}, a_{2}) \mu_{1}=k$ and
$\displaystyle \tau\mu_{1} \neq\frac{\pi}{2}+2n\pi$ for $n \geq0,$ we find that $\beta (a_{1},a_{2})<1.$ Then taking
into account \eqref{taucon} we obtain
\[
k<\lambda_{k,n}\mu_{1}=\frac{k\tau\mu_{1}}{\frac{\pi}{2}+2n\pi} < k
\frac{\frac{\pi}{2} + 2(n_{1} +1) \pi}{\frac{\pi}{2}+2n\pi}.
\]
Therefore we have $n\leq n_{1}.$ On the other hand, we have by the definition
that
\begin{equation}
\label{nierw}\frac{\tau}{2(j+1)\pi} \leq\lambda_{k,n} =\frac{k \tau}%
{(\frac{\pi}{2}+2n\pi)}\leq\frac{\tau}{2j\pi},
\end{equation}
which is equivalent to $kj \leq n < k(j+1).$

Therefore $1\leq k\leq[ n_{1}/j].$ Then noting that $\mu_{1}\neq\mu_{2},$ we
have that
\[
(a(t),\lambda)=((a_{1}(t),a_{2}(t)),\lambda)=((c_{1,k} \cos kt,0),\lambda
_{k,n} ) \quad\text{ for some } 1 \leq k \leq[ n_{1}/j], 1 \leq n \leq n_{1},
\]
where $c_{1,k}>0$ is such that $\beta(c_{1,k}\cos kt,0) \lambda_{k,n}\mu
_{1}=k. $

Similarly, we have that if $a_{2}\not \equiv 0,$%
\[
(a(t),\lambda)=((a_{1}(t),a_{2}(t)),\lambda)=((0,c_{2,k}\cos kt)),\lambda
_{k,n}) \quad\text{for some } 1 \leq k \leq[ n_{2}/j], 1 \leq n \leq n_{2},
\]
where where $c_{2,k}>0$ is such that $\beta(0,c_{2,k}\cos kt) \lambda_{k,n}
\mu_{2}=k.$

It is clear that the map $s \rightarrow\beta(su)$ is decreasing for any $u \in
D_{R}(E) \setminus cl(D_{r}(E)).$ Then since $\beta(a_{1}, a_{2}) <1,$ the map
$s \rightarrow\beta(s a_{1}, s a_{2})$ is strictly decreasing in
$[1-\varepsilon,1+\varepsilon].$

This implies that each $\left\{  (\rho(e^{i\theta},(a_{1}(t),a_{2}%
(t))),\lambda):\theta\in[0,2\pi)\right\}  $ is an isolated orbit satisfying
(\ref{Trans}).

Now fix $(a_{0}(t),\lambda_{k,n})=((a_{1}(t),a_{2}(t)),\lambda_{k,n}%
)=((c_{1,k}\cos kt,0),\lambda_{k,n}),$ where $1\leq k\leq[ n_{1}/j]$ and $1
\leq n < n_{1}.$ Then $(\dot a_{0}(t),0) =((\dot a_{1}(t),\dot a_{2}%
(t)),0)=((-c_{1,k} k \sin kt,0),0)$ is the tangent vector to the orbit $S^1 \cdot(a_{0},\lambda_{k,n})$ at
$(a_{0},\lambda_{k,n}).$

Summing up, we have proved that $(Q-G(\cdot,1))^{-1}(0)$ consists of a finite number of $S^1$-orbits. Below we prove
that these orbits are non-degenerate.

For simplicity of notation we put $x=(x_{1}(t),x_{2}(t)).$ Then
\begin{equation}
\label{fff2}f(x,\lambda)=\left(
\begin{array}
[c]{c}%
f_{1}(x,\lambda)\\
f_{2}(x,\lambda)
\end{array}
\right)  = \left(
\begin{array}
[c]{c}%
\displaystyle \int_{0}^{t} -\sigma(\lambda) \beta(x) \mu_{1} x_{1}(s -
\frac{\tau}{\lambda}) ds\\
\displaystyle \int_{0}^{t} -\sigma(\lambda) \beta(x) \mu_{2} x_{2}(s -
\frac{\tau}{\lambda}) ds
\end{array}
\right) ,
\end{equation}

\[
D_{x} f=\left(
\begin{array}
[c]{c}%
D_{x} f_{1}\\
D_{x} f_{2}%
\end{array}
\right) ,\: D_{\lambda}f=\left(
\begin{array}
[c]{c}%
D_{\lambda}f_{1}\\
D_{\lambda}f_{2}%
\end{array}
\right) .
\]

Then we obtain
\begin{equation}
\label{dxx2}D_{x}f(a_{0},\lambda_{k,n})(v)=
\end{equation}

\[
=- \int_{0}^{t} \left(
\begin{array}
[c]{c}%
\displaystyle k v_{1}(s- \frac{\tau}{\lambda_{k,n}})+ 2 \lambda_{k,n} \mu_{1}
a_{1}(s - \frac{\tau}{\lambda_{k,n}}) \xi^{\prime}(a_{0})\langle a_{0}, v
\rangle\\
\displaystyle k \frac{\mu_{2}}{\mu_{1}} v_{2}(s- \frac{\tau}{\lambda_{k,n}})
\end{array}
\right) \: ds
\]

\[
=- \int_{0}^{t} \left(
\begin{array}
[c]{c}%
\displaystyle k v_{1}(s- \frac{\tau}{\lambda_{k,n}})+ 2 \lambda_{k,n} \mu_{1}
a_{1}(s - \frac{\pi}{2k}) \xi^{\prime}(a_{0})\langle a_{1}, v_{1} \rangle\\
\displaystyle k \frac{\mu_{2}}{\mu_{1}} v_{2}(s- \frac{\tau}{\lambda_{k,n}})
\end{array}
\right) \: ds
\]

\[
=\left(
\begin{array}
[c]{cc}%
T_{1} & 0\\
0 & T_{2}%
\end{array}
\right) \left(
\begin{array}
[c]{c}%
v_{1}\\
v_{2}%
\end{array}
\right) ,
\]
where $\xi^{\prime}(a_{0}) < 0.$

\noindent On the other hand
\[
D_{\lambda}f(x,\lambda)=\left(
\begin{array}
[c]{c}%
\displaystyle \int_{0}^{t} - \beta(x) \mu_{1}(\sigma^{\prime}(\lambda) x_{1}(s
- \frac{\tau}{\lambda})+ \frac{\sigma(\lambda)\tau}{\lambda^{2}} \dot x_{1}(s
- \frac{\tau}{\lambda})) ds\\
\displaystyle \int_{0}^{t} - \beta(x) \mu_{2}(\sigma^{\prime}(\lambda) x_{2}(s
- \frac{\tau}{\lambda})+ \frac{\sigma(\lambda)\tau}{\lambda^{2}} \dot x_{2}(s
- \frac{\tau}{\lambda})) ds
\end{array}
\right) ,
\]
and noting that $\sigma(\lambda_{k,n})=\lambda_{k,n},\sigma^{\prime}%
(\lambda_{k,n})=0$ and that $a_{0}=(a_{1},0)$ we obtain
\begin{equation}
\label{dla2}D_{\lambda}f(a_{0},\lambda_{k,n})=\left(
\begin{array}
[c]{c}%
\displaystyle - \int_{0}^{t} \beta(a_{0}) \mu_{1} \lambda_{k,n}\frac{\tau
}{\lambda_{k,n}^{2}} \dot a_{1}(s-\frac{\tau}{\lambda_{k,n}}) ds\\
0
\end{array}
\right) =
\end{equation}
\[
= \left(
\begin{array}
[c]{c}%
\displaystyle - \int_{0}^{t} \frac{k\tau}{\lambda_{k,n}^{2}} \dot
a_{1}(s-\frac{\pi}{2k}) ds\\
0
\end{array}
\right) =\left(
\begin{array}
[c]{c}%
\displaystyle -\frac{k \tau}{\lambda_{k,n}^{2}} a_{1}(t - \frac{\pi}{2k})\\
0
\end{array}
\right)  = \left(
\begin{array}
[c]{c}%
T_{3}\\
0
\end{array}
\right) .
\]

Let us consider the following eigenvalue problem $v=\mu D_{x}f(a_{0}%
,\lambda_{k,n})v$ i.e.
\[
\left(
\begin{array}
[c]{c}%
v_{1}\\
v_{2}%
\end{array}
\right) v=\mu\left(
\begin{array}
[c]{c}%
T_{1} v_{1}\\
T_{2} v_{2}%
\end{array}
\right)
\]
for $v=(v_{1},v_{2}) \in E,$ which is equivalent to the following system%

\begin{equation}
\label{epr2}\left\{
\begin{array}
[c]{rl}%
\displaystyle \dot v_{1}(t)= & \displaystyle - \mu(\displaystyle k v_{1}(t-
\frac{\tau}{\lambda_{k,n}})+ 2 \lambda_{k,n} \mu_{1} a_{1}(t - \frac{\pi}{2k})
\xi^{\prime}(a_{0})\langle a_{1}, v_{1} \rangle),\\
\displaystyle \dot v_{2}(t)= & \displaystyle - \mu\frac{k\mu_{2}}{\mu_{1}}
v_{2}(t- \frac{\tau}{\lambda_{k,n}}).
\end{array}
\right.
\end{equation}

Since $\langle a_{1},\dot a_{1} \rangle=0,$ it is easy to verify that $\mu=1$
is the eigenvalue with corresponding eigenvector $(\dot a_{1},0).$

Summing up, we obtain
\[
Q - D f(a_{0},\lambda_{k,n})=Q - (D_{x} f(a_{0},\lambda_{k,n}),D_{\lambda}
f(a_{0},\lambda_{k,n}))=
\]
\[
=Q - T =\left(
\begin{array}
[c]{ccc}%
Id & 0 & 0\\
0 & Id & 0
\end{array}
\right)  - \left(
\begin{array}
[c]{ccc}%
T_{1} & 0 & T_{3}\\
0 & T_{2} & 0
\end{array}
\right)  : E \times\mathbb{R} \rightarrow E
\]
is a surjection such that $\ker(Q - T) = \mathrm{span} \{((\dot a_{1},0),0)\}$. Notice that we have just proved that
$S^1$-orbits of $(Q-G(\cdot,1))^{-1}(0)$ are non-degenerate.

In other words the assumptions of Theorem \ref{dgjm2} are satisfied.

Since $(a_{0}(t),\lambda_{k,n})$ is an isolated non-degenerate solution of \eqref{Trans} and
$S^1_{a_{0}}=\mathbb{Z}_{k},$  applying Theorem \ref{dgjm2} we obtain
\begin{equation}
\label{ghg}\mathrm{deg}_{\mathbb{Z}_{k}}(Q-f, \Omega)=\pm 1, \mathrm{deg}_{\mathbb{Z}_{k'}}(Q-f, \Omega)=0 \:
\mathrm{ for } \: k' > k
\end{equation}
for an open, bounded $S^1$-invariant subset $\Omega\subset cl(\Omega) \subset\Omega_{\lambda_{1}, \lambda_{2}}$ such
that $(Q-f)^{-1}(0) \cap\Omega=S^1 \cdot a_{0}\times\{\lambda_{k,n}\}.$

The same computation one can perform for

\noindent$(a_{0}(t),\lambda)=((a_{1}(t),a_{2}(t)),\lambda)=((0,c_{2,k}\cos kt)),\lambda_{k,n}) $  for some $1 \leq
k, n \leq n_{2},$ satisfying \eqref{nierw}.

Summing up, we have proved that  $(Q-f)^{-1}(0) \cap \Omega_{\lambda_1 \lambda_2}=S^1 \cdot a_1 \cup \ldots \cup S^1
\cdot a_p$ i.e. it consist of a finite number of non-degenerate $S^1$-orbits $S^1 \cdot a_1, \ldots, S^1  \cdot
a_p$, each with nontrivial $S^1$-index,   see formula \eqref{ghg} . Since these orbits are non-degenerate, there are
open bounded $S^1$-invariant subsets $U_i \subset cl(U_i) \subset \Omega_{\lambda_1 \lambda_2},i=1,\ldots,p,$ such
that $(Q -G(\cdot,1))^{-1}(0) \cap U_i=S^1  \cdot a_i, i=1, \ldots,p.$ And consequently by the properties of
$S^1$-degree we obtain
$$\mathrm{Deg}(Q-f,\Omega_{\lambda_{1},\lambda_{2}})=$$ $$= \mathrm{Deg}(Q-f,U_1)+ \ldots +
\mathrm{Deg}(Q-f,U_p) \in \Gamma.$$

Notice that for $\displaystyle k_{0}=\left[ \frac{n_{2}}{j} \right] $ only $n=n_{2}$ satisfies \eqref{nierw}.
Therefore there is exactly one solution of \eqref{Trans} in $\Omega_{\lambda_{1},\lambda_{2}}$ of the form
$((0,c_{2,k_{0} }\cos k_{0} t)),\lambda_{k_{0},n_{2}}).$ Moreover, other solutions of \eqref{Trans} are of the form
$((c_{1,k} \cos kt,0),\lambda _{k,n}) \text{ or } ((0,c_{2,k}\cos kt)),\lambda_{k,n}),$ where $k < k_{0}.$

In other words there is exactly one orbit with isotropy group $\mathbb{Z}_{k_0}$

Finally, combining Theorem \ref{dgjm2} with  \eqref{ghg} we obtain
$\mathrm{deg}_{\mathbb{Z}_{k_{0}}}(Q-f,\Omega_{\lambda_{1},\lambda_{2}}) \neq0,$ which completes the proof.
\end{proof}

\noindent\textbf{Proof of Theorem \ref{t}.}
 Without loss of generality we can assume that $n_1 < n_2.$ Fix $j \in \Phi(n_1,n_2)$ and define
$\lambda_1=\frac{\tau}{2(j+1)\pi},\lambda_2=\frac{\tau}{2j\pi},$ $\Omega_{\lambda_1 \lambda_2}=\Theta \times
(\lambda_1,\lambda_2),$ where $\Theta \subset E$ is an open bounded $S^1$-invariant subset defined by \eqref{Theta}.
In other words $\Omega_{\lambda_1 \lambda_2}$ is a cartesian product of an "annulus" $\Theta$ and an open interval
$(\lambda_1,\lambda_2).$

To complete the proof it is enough to show that $(Q -F(\cdot,1))^{-1}(0) \cap \Omega_{\lambda_1\lambda_2} \neq
\emptyset,$ where the operator $F$ is defined by formula \eqref{opef}. By Theorem \ref{dgjm} it is enough to show
that either $(Q -F(\cdot,1))^{-1}(0) \cap
\partial \Omega_{\lambda_1\lambda_2} \neq \emptyset,$ or
$\mathrm{Deg}(Q -F(\cdot,1),\Omega_{\lambda_{1},\lambda_{2}}) \neq \Theta \in \Gamma.$

By Lemma  \ref{homotopy1} we obtain that $(Q-F)^{-1}(0) \cap (\partial \Omega_{\lambda_1\lambda_2} \times [0,1]) =
\emptyset.$ Therefore by the homotopy property of the $S^1$-degree we obtain that $$\mathrm{Deg} (Q -F(\cdot,1),
\Omega_{\lambda_{1},\lambda_{2}}) = \mathrm{Deg}(Q -F(\cdot,0),\Omega_{\lambda_{1},\lambda_{2}}).$$ By Theorem
\ref{dgjm}, what is left is to show that $\mathrm{Deg}(Q-F(\cdot,0),\Omega_{\lambda_{1},\lambda_{2}}) \neq \Theta
\in \Gamma.$

From Lemma \ref{homotopy2} it follows that
$$\mathrm{Deg}(Q -F(\cdot,0),\Omega_{\lambda_{1},\lambda_{2}})=
\mathrm{Deg}(Q -G(\cdot,1),\Omega_{\lambda_{1},\lambda_{2}}).$$

Finally by Lemma \ref{homotopy3} we obtain $\mathrm{Deg}(Q -G(\cdot,1),\Omega_{\lambda_{1},\lambda_{2}}) \neq \Theta
\in \Gamma.$ Notice that we have just proved that  $\mathrm{Deg}(Q -F(\cdot,1),\Omega_{\lambda_{1},\lambda_{2}})
\neq \Theta \in \Gamma.$ The rest of the proof is a direct consequence of Theorem \ref{dgjm}.

\bibliographystyle{aplain}
\bibliography{acompat,mybibor11,mybibor11}

\begin{thebibliography}{9}                                                                                                %
\bibitem {[BAFAKR]}Z. Balanov, M. Farzamirad \& W. Krawcewicz,
\textit{Symmetric systems of van der Pol Equations}, Topol. Meth. Nonlin. Anal. \textbf{27(1)} (2005), 29-90,

\bibitem {[DGJM]}G. Dylawerski, K. Geba, J. Jodel \& W. Marzantowicz,
\textit{An $S^1$-equivariant degree and the Fuller index}, Ann. Pol. Math. \textbf{63} (1991), 243-280,

\bibitem {[GOPA]}K. Goparlsamy, \textit{Stability and oscillations in delay
differential equations of polulation dynamics}, Kluwer Academic Publishers, 1992,

\bibitem {[HALE1]}J. Hale, \textit{Theory of functional differential
equations}, Springer-Verlag, 1976,

\bibitem {[HR]}N. Hirano \& S. Rybicki, \textit{Existence of limit cycles for
coupled van der Pol Equations}, J. Diff. Equat. \textbf{195(1)} (2003), 194-209.
\end{thebibliography}

\end{document}